\documentclass[10pt]{amsart}
\usepackage{amssymb}
\usepackage{enumerate}
\usepackage[dvips]{graphicx}
\usepackage{epsf}
\usepackage{psfrag}
\usepackage{graphics}
\DeclareGraphicsExtensions{.eps,.art,.ART,.ps}
\newcommand{\bbR}{\mathbb{R}}      
\newcommand{\bbN}{\mathbb{N}}      
\newcommand{\bbZ}{\mathbb{Z}}      

\newcommand{\tr}{\operatorname{tr}}
\newcommand{\grad}{\operatorname{grad}}

\newcommand{\dive}{\operatorname{div}}
\newcommand{\diag}{\operatorname{diag}}
\newcommand{\cX}{\mathfrak{X}}
\newcommand{\cT}{\mathcal{T}}
\newtheorem{Thm}{Theorem}[section]
\newtheorem{Prop}[Thm]{Proposition}

\newtheorem{Cor}[Thm]{Corollary}
\theoremstyle{definition}
\newtheorem{Def}[Thm]{Definition}
\newtheorem{Example}[Thm]{Example}
\newtheorem{Remark}[Thm]{Remark}
\theoremstyle{remark}

\begin{document}
\begin{abstract}
We summarize the main ideas of General Relativity and Lorentzian geometry, leading to a proof of the simplest of the celebrated Hawking-Penrose singularity theorems. The reader is assumed to be familiar with Riemannian geometry and point set topology.
\end{abstract}
%
%
\title[Relativity and Singularities]{Relativity and Singularities -- A Short Introduction for Mathematicians}
\author{Jos\'{e} Nat\'{a}rio}
\address{Department of Mathematics, Instituto Superior T\'{e}cnico, Portugal}
\thanks{This work was partially supported by FCT/POCTI/FEDER}
\maketitle
%
%
%
\section*{Introduction}
Historically, much of the development of Riemannian geometry has been driven by General Relativity. This theory models spacetime as a Lorentzian manifold, which is analogous to a Riemannian manifold except that the positive definite metric is replaced by a metric with signature $(-,+, \ldots, +)$. Not only is Lorentzian geometry similar to Riemannian geometry in many respects but also Riemannian manifolds arise naturally as submanifolds of Lorentzian manifolds. Physical considerations then give rise to conjectures in Riemannian geometry. Recent examples of results inspired by such conjectures include the mass positivity theorem (Schoen and Yau, \cite{SY79}, \cite{SY81}) and the Penrose inequality (Bray, \cite{B01}, Huisken and Ilmanen, \cite{HI01}).

On the other hand, the effort involved in learning Lorentzian geometry is minimal once one has mastered Riemannian geometry. It therefore seems strange that many mathematicians (even geometers) choose not to do so. This may be in part due to the fact that most introductions to General Relativity start from first principles, developing the required differential geometry tools at length, and mostly focus on physical implications of the theory. A mathematician might prefer a shorter introduction to the subject from a more advanced starting point, focusing on interesting mathematical ideas. This paper aims to provide such an introduction, leading to a nontrivial result -- the simplest of the Hawking-Penrose singularity theorems (\cite{P65}, \cite{H67}, \cite{HP70}). These theorems basically state that physically reasonable Lorentzian manifolds (in a precise mathematical sense) must be singular (i.e.~geodesically incomplete). Since the motions of free-falling particles are represented by geodesics, this has the physical interpretation that General Relativity cannot be a complete description of Nature.

The paper is divided into three sections. The first section contains basic ideas of General Relativity and Lorentzian geometry: timelike, spacelike and null vectors and curves, matter models, the Einstein equation and its simplest solutions. Causality theory is developed in the second section, where we discuss time orientation, chronological and causal future and past sets, local causal structure, local maximizing properties of timelike geodesics (the Twin Paradox), the chronology condition, stable causality, domains of dependence and global hyperbolicity. The third section contains the proof of the singularity theorem. The proof has three ingredients: the first is that timelike geodesics cease to maximize the distance to a given time slice $S$ once a conjugate point is reached. The second is that, under a physically reasonable condition (the strong energy condition), conjugate points always occur. The third is that a length maximizing geodesic connecting $S$ to any given point $p$ always exists in a globally hyperbolic Lorentzian manifold. This is the most difficult (and mathematically interesting) point to prove; the proof is achieved by showing that the set of timelike curves connecting $S$ to $p$ with the Hausdorff metric is a compact space where the length functional is upper semicontinuous.

We assume the reader to be familiar with elementary Riemannian geometry (as in \cite{Carmo93}) and point set topology (as in \cite{Munkres00}). For the reader whose interest is aroused by this short introduction, there are many excellent texts on General Relativity, usually containing also the relevant differential and Lorentzian geometry. These range from introductory (\cite{Schutz02}) to more advanced (\cite{W84}) to encyclopedic (\cite{MTW73}). More mathematically oriented treatments can be found in \cite{BEE96}, \cite{ONeill83} (\cite{GHL04} also contains a brief glance at pseudo-Riemannian geometry). Causality and the singularity theorems are treated in greater detail in \cite{Penrose87}, \cite{HE95}, \cite{Naber88}.
%
%
%
\section{General Relativity and Lorentzian Geometry}
General Relativity is the physical theory of space, time and gravitation. It models spacetime (i.e.~the set of all physical events) as a $4$-dimensional Lorentzian manifold (spacetimes with different numbers of dimensions are also considered, e.g. in String Theory).

\begin{Def}
An $n$-dimensional {\bf pseudo-Riemannian manifold} is a pair $(M,g)$, where $M$ is an $n$-dimensional differentiable manifold and $g$ is a symmetric, nondegenerate $2$-tensor field on $M$ (called the {\bf metric}). A pseudo-Riemannian manifold is said to be {\bf Riemannian} if $g$ has signature $(+ \ldots +)$, and is said to be {\bf Lorentzian} if $g$ has signature $(-+ \ldots +)$.
\end{Def}

\begin{Example}
The simplest Riemannian manifold is {\bf Euclidean space}, which is $\bbR^n$ with the Riemannian metric
\[
g = dx^1 \otimes dx^1 + \ldots + dx^n \otimes dx^n.
\]
Analogously, the simplest Lorentzian manifold is {\bf Minkowski space}, which is $\bbR^{n+1}$ with the Lorentzian metric
\[
g = -dx^0 \otimes dx^0 + dx^1 \otimes dx^1 + \ldots + dx^n \otimes dx^n.
\]
\end{Example}

Lorentzian geometry is similar to Riemannian geometry in many respects. For instance, the Levi-Civita Theorem still holds (with the same proof) .

\begin{Thm}
{\bf (Levi-Civita)} Let $(M,g)$ be a pseudo-Riemannian manifold. Then there exists a unique symmetric affine connection compatible with the metric.
\end{Thm}

In particular, a Lorentzian manifold comes equipped with geodesics.

\begin{Example}
The Levi-Civita connection of Minkowski space is the trivial connection, and its geodesics are straight lines.
\end{Example}

On the other hand, the minus sign in the signature does introduce many novel features to Lorentzian geometry.

\begin{Def}
Let $(M,g)$ be a Lorentzian manifold and $p \in M$. A vector $v \in T_pM$ is said to be
\begin{enumerate}[(i)]
\item
{\bf Timelike} if $\langle v, v \rangle < 0$;
\item
{\bf Spacelike} if $\langle v, v \rangle > 0$;
\item
{\bf Null} if $\langle v, v \rangle = 0$.
\end{enumerate}
The {\bf length} of $v$ is $|v|=|\langle v, v \rangle|^\frac12$. (As usual one writes $\langle v, w \rangle$ for $g(v,w)$).
\end{Def}


A curve $c:I \subset \bbR \to M$ is said to be timelike, spacelike or null if its tangent vector $\dot{c}(t)$ is timelike, spacelike or null for all $t \in I$. If $c$ is a geodesic, then the nature of its tangent vector cannot change, as
\[
\frac{d}{dt} \langle \dot{c}, \dot{c} \rangle = 2 \left\langle \frac{D\dot{c}}{dt}, \dot{c} \right\rangle = 0.
\]

\begin{Example}
Interpreting the $x^0$-coordinate of Minkowski space as the time measured in some inertial frame, we see that timelike curves represent motions of particles such that
\[
\left(\frac{dx^1}{dx^0}\right)^2 + \ldots + \left(\frac{dx^n}{dx^0}\right)^2 < 1.
\]
It is assumed that units have been chosen so that $1$ is the maximum allowed velocity for a material particle (the speed of light). Therefore, timelike curves represent motions of material particles. Timelike geodesics, on the other hand, represent straight line motions with constant speed, i.e.~motions of {\bf free} particles. In addition, the length
\[
\tau(c)=\int_a^b \left| \dot{c} \right| dt
\]
of a timelike curve $c:[a,b] \to M$ is interpreted as the {\bf proper time} measured by the particle between events $c(a)$ and $c(b)$.

Null curves, on the other hand, represent motions at the speed of light, and null geodesics represent the motions of light rays.
\end{Example}

\begin{figure}[h!]
\begin{center}
\psfrag{timelike}{timelike}
\psfrag{null}{null}
\psfrag{spacelike}{spacelike}
\psfrag{timelike curve}{timelike curve}
\psfrag{spacelike geodesic}{spacelike geodesic}
\psfrag{t}{$x^0$}
\psfrag{x}{$x^1$}
\psfrag{y}{$x^n$}
\epsfxsize=.8\textwidth
\leavevmode
\epsfbox{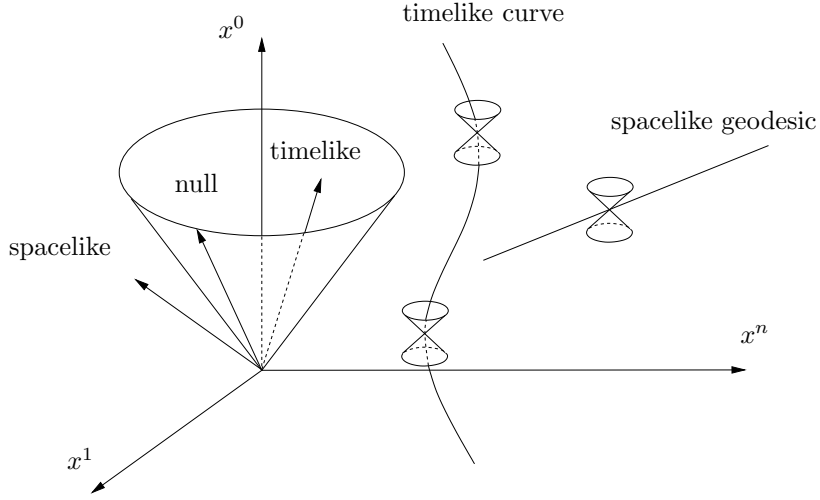}
\end{center}
\caption{Minkowski space.} 
\end{figure}

Einstein's great insight was realizing that to represent the gravitational field one must allow for general (curved) Lorentzian manifolds. The interpretations of timelike and null curves remain the same as in Minkowski space, except that timelike geodesics now represent {\bf free-falling} material particles.

This first step, however, is just the easy half of the problem. To complete the theory, one needs to decide which Lorentzian manifolds represent actual (physical) gravitational fields. It took Einstein several years of hard work to find the answer.

\begin{Def}
A Lorentzian manifold $(M,g)$ is said to satisfy the {\bf Einstein equation} if its Ricci curvature tensor satisfies
\[
Ric - \frac{S}2 g = T,
\]
where $S = \tr Ric$ is the scalar curvature and $T$ is the {\bf energy-momentum tensor} for a matter model in $(M,g)$.
\end{Def}

\begin{Example} \label{matter}
Rather than describing in detail what a matter model is, we simply list the simplest examples:
\begin{enumerate}[(i)]
\item
{\bf Vacuum:} Corresponds to taking $T=0$.
\item
{\bf Cosmological constant:} Corresponds to taking
\[
T = -\Lambda g
\]
for some constant $\Lambda \in \bbR$.
\item
{\bf Pressureless perfect fluid:} Is described by a {\bf rest mass density function} $\rho \in C^\infty(M)$ and a timelike unit {\bf velocity field} $U \in \cX(M)$ (whose integral lines are the motions of the fluid particles), and corresponds to the energy-momentum tensor
\[
T = \rho \, U^\sharp \otimes U^\sharp 
\]
($U^\sharp$ is the $1$-form associated to $U$ by the metric).
\end{enumerate}
\end{Example}

\begin{Example} \label{spacetimes}
We now list the simplest solutions $(M,g)$ of Einstein's field equation:
\begin{enumerate}[(i)]
\item
{\bf Minkowski space (1907):} Is clearly a vacuum solution, describing an universe without gravitation.
\item
{\bf Schwarzschild solution (1916):} Any vacuum solution admitting $O(n)$ as an isometry group is locally isometric to $M=\bbR^2 \times S^{n-1}$ with the metric
\[
g = - \left( 1 - \left(\frac{r_s}{r}\right)^{n-2} \right) dt \otimes dt + \left( 1 - \left(\frac{r_s}{r}\right)^{n-2} \right)^{-1} dr \otimes dr + r^2 h
\]
for some $r_s \in \bbR$, where $h$ is the round metric in $S^{n-1}$. If $r_s > 0$ then $r$ varies on either $(0,r_s)$ ({\bf Schwarzschild interior}) or $(r_s,+\infty)$ ({\bf Schwarzschild exterior}). These two regions can also be glued together by adding a {\bf horizon}, where $r=r_s$, in which case one obtains a model of a {\bf black hole}. Notice that in the Schwarzschild interior $r$ is a time coordinate, and the integral lines of $\frac{\partial}{\partial r}$ are {\bf incomplete} timelike geodesics, as as they cannot be continued past the curvature singularity at $r=0$.

For $n \geq 3$ this metric admits circular orbits, i.e.~geodesics of constant $r$ coordinate traversing a great circle in $S^{n-1}$, and these satisfy
\[
\left( \frac{d\varphi}{dt} \right)^2 = \frac{(n-2){r_s}^{n-2}}{2r^n}
\]
($\varphi$ being the angular coordinate in the great circle), which is exactly the condition for circular orbits around a point mass proportional to ${r_s}^{n-2}$ in Newtonian mechanics.
\item
{\bf de Sitter/Anti-de Sitter space (1917):} de Sitter space is simply the $(n+1)$-dimensional hyperboloid
\[
-(x^0)^2 + (x^1)^2 + \ldots + (x^{n+1})^2 = \alpha^2
\]
in $(n+2)$-dimensional Minkowski space, with the induced metric. Notice that this is the analogue of an Euclidean sphere, and is indeed a space of constant curvature. (As an aside, we remark that by reversing the sign of $\alpha^2$ one obtains two copies of the $n+1$-dimensional hyperbolic space of radius $\alpha$). de Sitter space is a solution of the Einstein equation with cosmological constant 
\[
\Lambda=\frac{n(n-1)}{2\alpha^2}.
\]
The induced metric can be written as
\[
g = \alpha^2 \left( - dt \otimes dt + \cosh^2 t \,\, h \right),
\]
where $x^0=\alpha\sinh t$ and $h$ is the round metric in $S^n$. Therefore one can think of de Sitter space as a spherical universe which contracts to a minimum radius $\alpha$ and then re-expands at an exponential rate. This cosmic repulsion is driven by the positive cosmological constant.

Anti-de Sitter space is the universal cover of the $(n+1)$-dimensional hyperboloid
\[
(x^1)^2 + \ldots + (x^n)^2 - (x^{n+1})^2 - (x^{n+2})^2 = - \alpha^2
\]
in $\bbR^{n+2}$ with the pseudo-Riemannian metric
\[
dx^1 \otimes dx^1 + \ldots + dx^n \otimes dx^n - dx^{n+1} \otimes dx^{n+1} - dx^{n+2} \otimes dx^{n+2},
\]
and is a solution of the Einstein equation with cosmological constant 
\[
\Lambda=-\frac{n(n-1)}{2\alpha^2}
\]
(again being a space of constant curvature). The induced metric can be written as
\begin{align*}
g & =  \alpha^2 \left( - \cosh^2 \xi \, dt \otimes dt +  d\xi \otimes d\xi +\sinh^2 \xi  \,\, h \right) \\
& = \frac{\alpha^2}{\cos^2 x} \left(- dt \otimes dt + dx \otimes dx + \sin^2 x \,\, h \right)
\end{align*}
where $(t,\xi) \in \bbR \times \left(0,+\infty\right)$ are defined by 
\[
\begin{cases}
x^{n+1}=\alpha \cosh \xi \cos t\\
x^{n+2}=\alpha \cosh \xi \sin t
\end{cases},
\]
$x \in \left(0,\frac\pi2\right)$ satisfies $\cos x = \frac1{\cosh \xi}$, and $h$ is the round metric in $S^{n-1}$. Therefore one can think of de Sitter space as a static universe whose spatial sections are hyperbolic spaces (hence conformal to half spheres).

\item
{\bf Friedmann-Lema\^\i tre-Robertson-Walker models (1922):} If $(\Sigma,h)$ is a $n$-dimensional manifold of constant curvature $k=-1,0,1$ (i.e.~hyperbolic space, Euclidean space, $S^n$ or quotients of these by discrete groups of isometries) then $(\bbR \times \Sigma, -dt \otimes dt + a^2(t) h)$ is a solution of the Einstein equation (representing an expanding or contracting universe) for a perfect fluid with velocity field $U = \frac{\partial}{\partial t}$ and rest mass density
\[
\rho = \frac{n(n-1) \alpha}{a^n},
\]
where the function $a(t)$ satisfies the first order ODE
\[
\frac{\dot{a}^2}{2} - \frac{\alpha}{a^{n-2}} = - \frac{k}2
\]
(for any constant $\alpha \in \bbR$). Notice that this is the equation of conservation of energy for a particle in the potential
\[
V(a) = - \frac{\alpha}{a^{n-2}}.
\]
Therefore, if $\alpha>0 \Leftrightarrow \rho > 0$ then $a(t)$ blows up in finite time, corresponding to a curvature singularity ({\bf Big Bang} or {\bf Big Crunch}). In this case, the timelike geodesics given by the integral curves of $\frac{\partial}{\partial t}$ are {\bf incomplete}.
\end{enumerate}
\end{Example}

\begin{Remark}
Let $(M,g)$ be a $(n+1)$-dimensional Lorentzian manifold. For $n=1$ one has the identity
\[
Ric - \frac{S}2 g = 0,
\]
and hence {\bf any} $2$-dimensional Lorentzian manifold is a vacuum solution of the Einstein equation (this being the only allowed matter model).

As is the case with all the fundamental equations of Physics, the Einstein equation can be derived from a variational principle. For instance, Lorentzian manifolds satisfying the vacuum Einstein equation are critical points of the {\bf Einstein-Hilbert action}
\[
\mathcal{A}[g] = \int_M S.
\]
The fact that all $2$-dimensional Lorentzian manifolds are vacuum solutions of the Einstein equation can thus be turned into a proof of the Gauss-Bonnet Theorem.

If $n \geq 2$, the trace of the Einstein equation  yields
\[
S = - \frac{2 \tr T}{n-1},
\]
and hence the Einstein equation can be rewritten as
\[
Ric = T - \frac{\tr T}{n-1} g.
\]
In particular, any solution of the Einstein equation must be Ricci-flat at points where $T=0$. This means that there is no gravitational field (i.e.~curvature) in vacuum for $n=2$.
\end{Remark}
%
%
\section{Causality}
We will now discuss the causal features of Lorentzian manifolds, a subject which has no parallel in Riemannian geometry. 

\begin{Def}
Let $(M,g)$ be a Lorentzian manifold and $p \in M$. Two timelike vectors $v,w \in T_pM$ are said to have the same (resp. opposite) {\bf time orientation} if $\langle v, w \rangle < 0$ (resp. $\langle v, w \rangle > 0$).
\end{Def}

Notice that $\langle v, w \rangle = 0$ cannot occur for timelike vectors.

\begin{Def}
A Lorentzian manifold $(M,g)$ is said to be {\bf time-orientable} if there exists a nonvanishing timelike vector field $T\in\cX(M)$.
\end{Def}

This means that one can consistently time orient all tangent spaces $T_pM$. If $(M,g)$ is connected and time-orientable then it has exactly two time orientations: nonvanishing timelike vector fields $T,U \in\cX(M)$ define the same time orientation if and only if $\langle T, U \rangle < 0$. Timelike vectors $v \in T_pM$ with the same time orientation as $T_p$ are said to be {\bf future-pointing} for the time orientation determined by $T$. This easily extends to nonvanishing null vectors.

\begin{Example}
The M\"obius band $\bbR \times [0,1] / \sim$, where $\sim$ is the equivalence relation $(t,0) \sim (-t, 1)$, admits the non-time-orientable Lorentzian metric
\[
g = - dt \otimes dt + dx \otimes dx,
\]
as well as the time-orientable Lorentzian metric $h=-g$.
\end{Example}

The usual proof for the existence of an orientable double cover for a non-orientable manifold can be easily adapted to prove

\begin{Prop}
Any non-time-orientable Lorentzian manifold $(M,g)$ has a {\bf time-orientable double cover}, i.e.~a time-orientable Lorentzian manifold $(\overline{M},\overline{g})$ and a local isometry $\pi:\overline{M}\to M$ such that every point in $M$ has two preimages by $\pi$.
\end{Prop}

A connected time-orientable Lorentzian manifold admits a nonvanishing vector field, and hence is either noncompact or has zero Euler characteristic. The same is true for a non-time-orientable Lorentzian manifold, for it must be true for its time-orientable double cover. On the other hand, if a connected differentiable manifold $M$ is either noncompact or has zero Euler characteristic then it admits a nonvanishing vector field $T$, which we can assume to be of unit length for some Riemannian metric $h$ on $M$. It is then easy to check that  $g = -2 T^\sharp \otimes T^\sharp + h$ is a Lorentzian metric on $M$. In other words:

\begin{Prop}
A connected differentiable manifold $M$ admits a Lorentzian metric if and only if is noncompact or has zero Euler characteristic.
\end{Prop}

\begin{Def}
Let $(M,g)$ be a time-oriented Lorentzian manifold.
\begin{enumerate}[(i)]
\item
A timelike curve $c:I \subset \bbR \to M$ is said to be {\bf future-directed} if $\dot{c}$ is future-pointing.
\item
The {\bf chronological future} of $p\in M$ is the set $I^+(p)$ of all points to which $p$ can be connected by a future-directed timelike curve.
\item
A {\bf future-directed causal curve} is a curve $c:I \subset \bbR \to M$ such that $\dot{c}$ is non-spacelike and future-pointing (if nonzero).
\item
The {\bf causal future} of $p\in M$ is the set $J^+(p)$ of all points to which $p$ can be connected by a future-directed causal curve.
\end{enumerate}
\end{Def}

One can make analogous definitions replacing ``future'' with ``past''. In general, the chronological and causal futures and pasts can be quite complicated sets, because of global features. Locally, however, causal properties are similar to those of Minkowski space. More precisely, we have the following statement:

\begin{Prop} \label{prop7.7.1}
Let $(M,g)$ be a time-oriented Lorentzian manifold. Then each point $p_0\in M$ has an open neighborhood $V \subset M$ such that the Lorentzian manifold $(V,g)$ satisfies:
\begin{enumerate}
\item If $p,q \in V$ then there exists a unique geodesic joining $p$ to $q$ (i.e.~$V$ is {\bf geodesically convex});
\item $q\in I^+(p)$ iff there exists a future-directed timelike geodesic connecting $p$ to $q$;
\item $J^+(p) = \overline{I^+(p)}$;
\item $q \in J^+(p)$ iff there exists a future-directed timelike or null geodesic connecting $p$ to $q$.
\end{enumerate}
\end{Prop}

\begin{proof}
The existence of geodesically convex neighborhoods holds for any affine connection (cf.~\cite{KN96}). Moreover, one can assume such neighborhoods to be {\bf totally normal}, i.e.~normal neighborhoods of all of their points.

To prove assertion (2), we start by noticing that if there exists a future-directed timelike geodesic connecting $p$ to $q$ then it is obvious that $q \in I^+(p)$. Suppose now that $q \in I^+(p)$; then there exists a future-directed timelike curve $c:[0,1] \to V$ such that $c(0)=p$ and $c(1)=q$. Choose normal coordinates $(x^0, x^1, \ldots, x^n)$ given by the parametrization
\[
\varphi(x^0,x^1,\ldots,x^n)=\exp_p(x^0 E_0 + x^1 E_1 + \ldots + x^n E_n),
\]
where $\{E_0,E_1,\ldots,E_n\}$ is an orthonormal basis of $T_pM$ (with $E_0$ timelike and future-pointing). These are global coordinates in $V$, since $V$ is totally normal. Defining
\begin{align*}
W_p(q) & := - (x^0(q))^2 + (x^1(q))^2 + \ldots + (x^n(q))^2 \\
& = \sum_{\mu,\nu=0}^n \eta_{\mu\nu}x^\mu(q)x^\nu(q),
\end{align*}
with $(\eta_{\mu\nu}) = \diag(-1,1,\ldots,1)$, we have to show that $W_p(q)<0$. Let $W_p(t):= W_p(c(t))$. Since $x^\mu(p)=0$ $(\mu = 0,1,\ldots,n)$, we have $W_p(0)=0$. Setting $x^\mu(t)=x^\mu(c(t))$, we have
\begin{align*}
& \dot{W}_p(t) = 2 \sum_{\mu,\nu=0}^n \eta_{\mu\nu}x^\mu(t)\dot{x}^\nu(t);\\
& \ddot{W}_p(t) = 2 \sum_{\mu,\nu=0}^n \eta_{\mu\nu}x^\mu(t)\ddot{x}^\nu(t) + 2\sum_{\mu,\nu=0}^n \eta_{\mu\nu}\dot{x}^\mu(t)\dot{x}^\nu(t),
\end{align*}
and consequently (recalling that $\left(d \exp_p\right)_p = \text{id}$)
\begin{align*}
& \dot{W}_p(0) = 0;\\
& \ddot{W}_p(0) = 2\langle \dot{c}(0), \dot{c}(0) \rangle < 0.
\end{align*}
Therefore there exists $\varepsilon > 0$ such that $W_p(t) < 0$ for $t \in (0, \varepsilon)$.

The same proof as in Riemannian geometry shows that the unit tangent vector field to timelike geodesics through $p$ is
\[
X = - \grad\left(-W_p\right)^\frac12 = \frac12 \left(-W_p\right)^{-\frac12} \grad W_p
\]
({\bf Gauss Lemma}), where the gradient of a function is defined as in Riemannian geometry (notice however that in Lorentzian geometry a smooth function $f$ {\bf decreases} along the direction of $\grad f$ if $\grad f$ is timelike). Consequently $\grad W_p$ is tangent to timelike geodesics through $p$, being future-pointing on future-directed timelike geodesics.

Suppose that $W_p(t) < 0$. Then 
\[
\dot{W}(t) = \left\langle \left(\grad W_p\right)_{c(t)}, \dot{c}(t) \right\rangle < 0
\]
as both $\left(\grad W_p\right)_{c(t)}$ and $\dot{c}(t)$ are timelike future-pointing. We conclude that we must have $W_p(t)<0$ for all $t\in[0,1]$. In particular, $W_p(q)=W_p(1)<0$, and hence there exists a future-directed timelike geodesic connecting $p$ to $q$.

Assertion (3) can be proved by using the global normal coordinates $(x^0,x^1,\ldots,x^n)$ of $V$ to approximate causal curves by timelike curves. Once this is done, (4) is obvious from the fact that $\exp_p$ is a diffeomorphism onto $V$.
\end{proof}

A simple application of Proposition~\ref{prop7.7.1} is proving

\begin{Prop}
Let $(M,g)$ be a time-oriented Lorentzian manifold and $p \in M$. Then:
\begin{enumerate}[(i)]
\item
If $q \in I^+(p)$ and $r \in I^+(q)$ then $r \in I^+(p)$;
\item
If $q \in J^+(p)$ and $r \in J^+(q)$ then $r \in J^+(p)$;
\item
$I^+(p)$ is an open set.
\end{enumerate}
\end{Prop}

In Lorentzian geometry there are no curves of minimal length, as any two points in the same connected component can be connected by piecewise null curves. However, there do exist curves with {\bf maximal} length, and these are timelike geodesics. More precisely, we have the following statement:

\begin{Prop} \label{prop7.7.2}
{\bf (Twin Paradox)} Let $(M,g)$ be a time-oriented Lorentzian manifold and $p_0 \in M$. Then there exists a geodesically convex open neighborhood $V \subset M$ of $p_0$ such that the Lorentzian manifold $(V,g)$ satisfies the following property: if $q \in I^+(p)$, $c$ is the timelike geodesic connecting $p$ to $q$ and $\gamma$ is any timelike curve connecting $p$ to $q$, then $\tau(\gamma) \leq \tau(c)$, with equality iff $\gamma$ is a  is a reparametrization of $c$.
\end{Prop}

\begin{proof}
Choose $V$ as in the proof of Proposition~\ref{prop7.7.1}. Any timelike curve $\gamma:[0,1] \to V$ satisfying $\gamma(0)=p$, $\gamma(1)=q$ can be written as
\[
\gamma(t)=\exp_p(r(t)n(t)),
\]
for $t \in [0,1]$, where $r(t) \geq 0$ and $\langle n(t), n(t) \rangle = -1$. We have
\[
\dot{\gamma}(t)=(\exp_p)_*\left(\dot{r}(t)n(t)+r(t)\dot{n}(t)\right).
\]
Since $\langle n(t), n(t) \rangle = -1$, we have $\langle \dot{n}(t), n(t) \rangle = 0$, and consequently $\dot{n}(t)$ is tangent to the level surfaces of the function $v \mapsto \langle v, v \rangle$. We conclude that
\[
\dot{\gamma}(t) = \dot{r}(t) X_{\gamma(t)} + Y(t),
\]
where $X$ is the unit tangent vector field to timelike geodesics through $p$ and $Y(t)=r(t)(\exp_p)_*\dot{n}(t)$ is tangent to the level surfaces of $W_p$ (hence orthogonal to $X_{\gamma(t)}$). Consequently,
\begin{align*}
\tau(\gamma) & = \int_0^1 \left|\left\langle \dot{r}(t) X_{\gamma(t)} + Y(t),\dot{r}(t) X_{\gamma(t)} + Y(t) \right\rangle\right|^\frac12 dt \\
& = \int_0^1 \left( \dot{r}(t)^2 - |Y(t)|^2 \right)^\frac12 dt \\
& \leq \int_0^1 \dot{r}(t) dt = r(1) = \tau(c),
\end{align*}
(where we've used the facts that $\dot{r}(t)> 0$ for all $t \in [0,1]$, as $\dot{c}$ is future-pointing, and $\tau(c)=r(1)$, as $q=\exp_p(r(1)n(1))$. It should be clear that $\tau(\gamma)=\tau(c)$ if and only if $|Y(t)|\equiv 0 \Leftrightarrow Y(t)\equiv 0$ ($Y(t)$ is spacelike) for all $t \in [0,1]$, implying that $n$ is constant. In this case, $\gamma(t)=\exp_p({r(t)n})$ is, up to reparametrization, the geodesic through $p$ with initial condition $n \in T_p M$.
\end{proof}

\begin{Remark}
Proposition~\ref{prop7.7.2} can be interpreted as follows: if two observers (e.g. two twins) meet at some event, are separated, and meet again at a later event, then the free-falling twin will always measure more time to have passed between the meetings.
\end{Remark}

For physical applications, it is important to demand that a Lorentzian manifold satisfies reasonable causality conditions. The simplest of these conditions excludes time travel, i.e.~the possibility of a particle returning to an event in its past history.

\begin{Def}
A Lorentzian manifold $(M,g)$ is said to satisfy the {\bf chronology condition} if it does not contain closed timelike curves.
\end{Def}

This condition is violated by compact Lorentzian manifolds:

\begin{Prop}
Any compact Lorentzian manifold $(M,g)$ contains closed timelike curves.
\end{Prop}

\begin{proof}
Taking if necessary the time-orientable double cover, we can assume that $(M,g)$ is time-oriented. It is easy to check that $\{ I^+(p) \}_{p \in M}$ is an open cover of $M$. If $M$ is compact, we can obtain a finite subcover $\{ I^+(p_1), \ldots, I^+(p_N) \}$. Now if $p_1 \in I^+(p_i)$ for $i \neq 1$ then $I^+(p_1) \subset I^+(p_i)$, and we can exclude $I^+(p_1)$ from the subcover. Therefore, we can assume without loss of generality that $p_1 \in I^+(p_1)$, and hence there exists a closed timelike curve starting and ending at $p_1$.
\end{proof}

A stronger restriction on the causal behavior of the Lorentzian manifold is the following:

\begin{Def}
A Lorentzian manifold $(M,g)$ is said to be {\bf stably causal} if there exists a {\bf global time function}, i.e.~a smooth function $t:M \to \bbR$ such that $\grad(t)$ is timelike.
\end{Def}

In particular, a stably causal Lorentzian manifold is time-orientable. We choose the time orientation defined by $-\grad(t)$, so that $t$ increases along future-directed timelike curves. Notice that this implies that no closed timelike curves can exist, i.e.~any stably causal Lorentzian manifold satisfies the chronology condition. In fact, any small perturbation of a causally stable Lorentzian manifold still satisfies the chronology condition.

\begin{Def}
Let $(M,g)$ be a time-oriented Lorentzian manifold. 
\begin{enumerate}[(i)]
\item
A smooth future-directed causal curve $c:(a,b) \to M$ (with possibly $a=-\infty$ or $b=+\infty$) is said to be {\bf future-inextendible} if $\lim_{t \to b} c(t)$ does not exist.
\item
The {\bf past domain of dependence} of $S\subset M$ is the set $D^-(S)$ of all points $p \in M$ such that any future-inextendible causal curve starting at $p$ intersects $S$.
\end{enumerate}
\end{Def}

One can make analogous definitions replacing ``future'' with ``past''. The {\bf domain of dependence} of $S$ is simply the set $D(S)=D^+(S)\cup D^-(S)$.

\begin{Def}
A Lorentzian manifold $(M,g)$ is said to be {\bf globally hyperbolic} if it is stably causal and there exists a time function $t:M \to \bbR$ such that the time slices $S_a = t^{-1}(a)$ are {\bf Cauchy hypersurfaces}, i.e.~satisfy $D(S_a)=M$.
\end{Def}

\begin{figure}[h!]
\begin{center}
\psfrag{t}{$t$}
\psfrag{x}{$x$}
\psfrag{S}{$S$}
\psfrag{D}{$D(S)$}
\psfrag{(-1,0)}{$(-1,0)$}
\psfrag{J+}{$J^+(-1,0)$}
\epsfxsize=1.0\textwidth
\leavevmode
\epsfbox{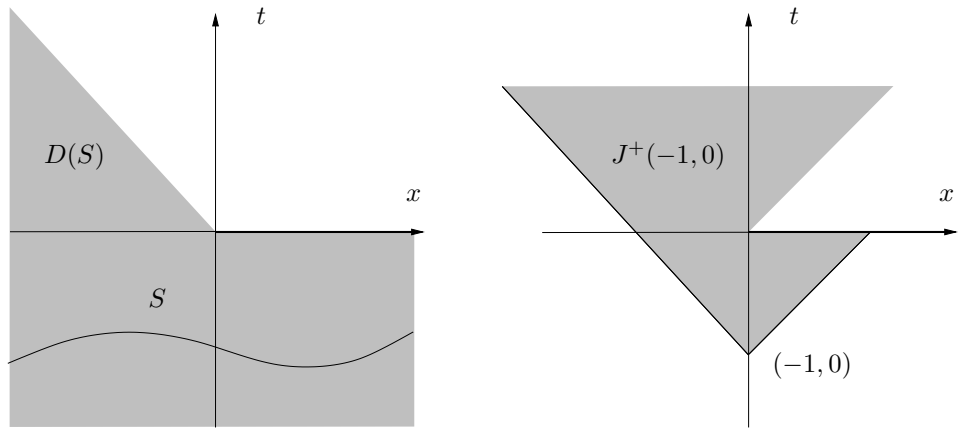}
\end{center}
\caption{The open set $U$ is not globally hyperbolic, and the causal future of $(-1,0)$ is not closed.} \label{Figure_not_gh}
\end{figure}

\begin{Example} \label{not_global_hyp}
The open set
\[
U = \{ (t,x) \in \bbR^2 \mid t \neq 0 \text{ or } x < 0 \}
\]
with the Minkowski metric
\[
g = - dt \otimes dt + dx \otimes dx
\]
is stably causal Lorentzian manifold (the coordinate $t:U \to \bbR$ is a global time function) which is not globally hyperbolic (cf.~Figure~\ref{Figure_not_gh}). Notice that $J^+(-1,0)$ is not closed. Moreover, the supremum of the lengths of timelike curves connecting $(-1,0)$ to $(1,0)$ is clearly $2$, but it is not attained by any timelike curve.
\end{Example}

%
%
%
\section{Singularity Theorem}
As we have seen in Example~\ref{spacetimes}, both the Schwarzschild solution and the Friedmann-Lema\^\i tre-Robertson-Walker cosmological models display singularities, beyond which timelike geodesics cannot be continued.

\begin{Def}
A Lorentzian manifold $(M,g)$ is {\bf singular} if it is not geodesically complete.
\end{Def}

One could think that the examples above are singular due to their high degree of symmetry, and that more realistic Lorentzian manifolds would be generically non-singular. We will show that this is not the case: any sufficiently small perturbation of these solutions will still be singular.

The question of whether a given Riemannian manifold is geodesically complete is settled by the Hopf-Rinow Theorem. Unfortunately, this theorem does not hold in Lorentzian geometry (essentially because one cannot use the metric to define a distance function). 

\begin{figure}[h!]
\begin{center}
\psfrag{t}{$t$}
\psfrag{x}{$x$}
\psfrag{(-p,0)}{$(-\pi,0)$}
\psfrag{(p,0)}{$(\pi,0)$}
\psfrag{p2}{$\frac{\pi}2$}
\psfrag{-p2}{$-\frac{\pi}2$}
\epsfxsize=.5\textwidth
\leavevmode
\epsfbox{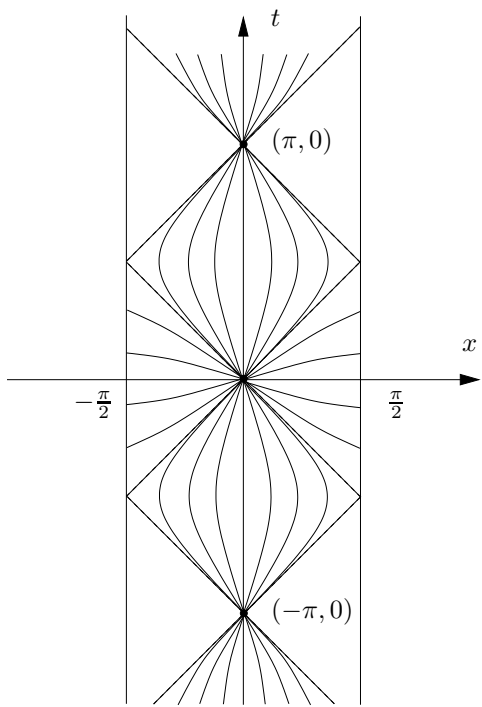}
\end{center}
\caption{The exponential map is not surjective in $2$-dimensional Anti-de Sitter space.} \label{AdS}
\end{figure}

\begin{Example} \hspace{1cm} \label{Hopf-Rinow}
\begin{enumerate}[(i)]
\item {\bf Clifton-Pohl torus:} Consider the Lorentzian metric
\[
\overline{g} = \frac1{u^2 + v^2} (du \otimes dv + dv \otimes du)
\]
on $\overline{M}=\bbR^2 \setminus \{ 0 \}$. The Lie group $\bbZ$ acts freely and properly on $\overline{M}$ by isometries through 
\[
n \cdot (u,v) = (2^n u, 2^n v),
\]
and this determines a Lorentzian metric $g$ on $M = \overline{M} / \bbZ \cong T^2$. A trivial calculation shows that there exist null geodesics satisfying $u \equiv 0$ and
\[
\frac{d}{dt}\left( \frac1{v^2} \frac{dv}{dt}\right) = 0.
\]
Since most solutions of this ODE blow up in finite time, we see that $(M,g)$ is not geodesically complete (although $M$ is compact).
\item \label{2dAdS} {\bf $2$-dimensional Anti-de Sitter space:} In spite of not solving the Einstein equation in $2$ dimensions, we can still consider the unit $2$-dimensional Anti-de Sitter space $(M,g)$, i.e.~the universal cover of the submanifold of $\bbR^3$
\[
\overline{M} = \{(u,v,w)\in \bbR^3 \mid u^2+v^2-w^2 = 1 )\}
\]
with the Lorentzian metric $\overline{g}$ induced by the pseudo-Riemannian metric
\[
- du\otimes du - dv\otimes dv + dw\otimes dw.
\]
As was seen in Example~\ref{spacetimes}, there exist global coordinates $(t,x) \in \bbR \times (-\frac\pi2, \frac\pi2)$ in $M$ such that and 
\[
g = \frac1{\cos^2 x} (- dt \otimes dt + dx \otimes dx)
\]
(notice the different range of the coordinate $x$, arising from the fact that $S^0=\{-1,1\}$). Therefore $(M,g)$ is conformal to the open subset $\bbR \times (-\frac\pi2, \frac\pi2)$ of Minkowski space, and its causality properties are the same. In particular, $(M,g)$ is not globally hyperbolic.

As is the case with the Euclidean $2$-sphere, the geodesics of $(\overline{M},\overline{g})$ can be obtained by intersecting $\overline{M}$ with $2$-planes through the origin. This fact can be used to prove that $(\overline{M},\overline{g})$, and hence $(M,g)$, are geodesically complete, and also that all timelike geodesics through the point $p$ with coordinates $(0,0)$ on $(M,g)$ refocus on the points with coordinates $(\pm \pi,0)$ (cf.~Figure~\ref{AdS}). On the other hand, spacelike geodesics through $p$ are timelike geodesics of the Lorentzian metric $-g$, and hence are confined to the chronological future and past of $p$ in this metric. Therefore $\exp_p$ is not surjective, although $(M,g)$ is geodesically complete.

Incidentally, if $\varepsilon > 0$ and $q$ is the point with coordinates $(\pi + \varepsilon, 0)$ then there exist piecewise smooth causal curves connecting $p$ to $q$ with arbitrarily large length: simply take a future-directed null geodesic from $p$ to the line $x=x_0$, a past-directed null geodesic from $q$ the same line, and the portion of this line between the two geodesics. The resulting curve has length greater than $\frac{\varepsilon}{\cos(x_0)}$, and this can be made arbitrarily large by making $x_0 \to \frac\pi2$. These curves can be easily smoothed into timelike curves with arbitrarily large length.
\end{enumerate}
\end{Example}

We now proceed to show that geodesic incompleteness is a generic feature of Lorentzian manifolds satisfying the strong energy condition.

\begin{Def}
A Lorentzian manifold $(M,g)$ is said to satisfy the {\bf strong energy condition} if $Ric(V,V) \geq 0$ for any timelike vector field $V \in \cX(M)$.
\end{Def}

If $n \geq 2$, this condition is equivalent (by the Einstein equation) to requiring that the energy-momentum tensor $T$ satisfies 
\[
T(V,V) \geq \frac{\tr T}{n-1} \langle V, V \rangle
\]
for any timelike vector field $V \in \cX(M)$, which turns out to be physically reasonable.

\begin{Example} Let us see the meaning of the strong energy condition for each of the matter models in Example~\ref{matter}:
\begin{enumerate}[(i)]
\item
{\bf Vacuum:} Trivially satisfied.
\item
{\bf Cosmological constant:} Equivalent to $\Lambda \leq 0$.
\item
{\bf Pressureless perfect fluid:} Becomes
\[
\rho \left( \left\langle U,V \right\rangle^2 + \frac1{n-1} \left\langle V,V \right\rangle \right) \geq 0,
\]
or, since the term in brackets is easily seen to be non-negative, simply $\rho \geq 0$.
\end{enumerate}
\end{Example}

\begin{Def}
Let $(M,g)$ be a globally hyperbolic Lorentzian manifold and $S$ a Cauchy hypersurface with future-pointing normal vector field $n$. Let $c_p$ be the timelike geodesic with initial condition $n_p$ for each point $p \in S$. The {\bf exponential map} $\exp:U \to M$ (defined on an open neighborhood $U$ of $\{0\} \times S \subset \bbR \times S$) is the map $\exp(t,p)=c_p(t)$. The critical values of $\exp$ are said to be {\bf conjugate points} to $S$.
\end{Def}

If the geodesic $c_p$ has no conjugate points between $c_p(0)=p$ and $c_p(t_0)$ then there exists an open neighborhood $V$ of $c_p([0,t_0])$ which can be foliated by images of geodesics orthogonal to $S$. The tangent vectors to these geodesics yield a unit timelike vector field $X \in \cX(V)$, which by the Gauss Lemma satisfies $X = - \grad t$, where $t:V \to \bbR$ is the distance along the geodesics. Therefore $X^\sharp = -dt$, and the covariant derivative $K=\nabla X^\sharp$ is a symmetric tensor, which by the geodesic equation must satisfy $\iota_X \nabla X^\sharp = \nabla_X X^\sharp = 0$.

\begin{Def}
The divergence $\theta = \dive X = \tr K$ of the vector field $X \in \cX(V)$ is called the {\bf expansion} of the family of timelike geodesics on $V$.
\end{Def}

\begin{Remark}
$K$ is actually the second fundamental form of the family of hypersurfaces of constant $t$, and $\theta$ gives the logarithmic rate of change of their volume element along $X$ (as is easily seen from the Divergence Theorem). 
\end{Remark}

\begin{Prop} \label{exist_conj_prop}
Let $(M,g)$ be a globally hyperbolic Lorentzian manifold satisfying the strong energy condition, $S$ a Cauchy hypersurface and $p \in S$ a point where $\theta=\theta_0 < 0$. If the geodesic $c_p$ can be extended to a distance $t_0= -\frac{n}{\theta_0}$ to the future of $S$, then it contains at least a point conjugate to $S$.
\end{Prop}

\begin{proof}
Suppose that there were no points in $c_p([0,t_0])$ conjugate to $S$. Then there would exist an open neighborhood $V$ of $c_p([0,t_0])$ as above. An easy calculation shows that the {\bf Raychaudhuri equation}
\[
X(\theta) + \tr K^2 + Ric(X,X) = 0
\]
holds. Since $(M,g)$ satisfies the strong energy condition, we have $R(X,X) \geq 0$, and hence
\[
X(\theta) + \tr K^2 \leq 0
\]
The Cauchy-Schwarz inequality implies that
\[
(\tr A)^2 \leq n \tr(A^tA),
\]
for any $n \times n$ matrix $A$. Since $K$ is symmetric and vanishes on $X$, we have
\[
\tr K^2 \geq \frac1n \theta^2.
\]
Setting $\theta(t):=\theta(c_p(t))$, we see that
\[
\frac{d \theta}{d t} + \frac1n \theta^2 \leq 0.
\] 
Integrating this inequality one obtains
\[
\frac1{\theta} \geq \frac1{\theta_0} + \frac{t}n,
\]
and hence $\theta$ must blow up at a value of $t$ no greater than $t_0=-\frac{n}{\theta_0}$. This yields a contradiction, as $\theta$ is smooth function on $V$.
\end{proof}

\begin{Prop} \label{conjugate_prop}
Let $(M,g)$ be a globally hyperbolic Lorentzian manifold, $S$ a Cauchy hypersurface, $p \in M$ and $c$ a timelike geodesic through $p$ orthogonal to $S$. If there exists a conjugate point between $S$ and $p$ then $c$ does not maximize length (among the timelike curves connecting $S$ to $p$).
\end{Prop}

\begin{figure}[h!]
\begin{center}
\psfrag{p}{$p$}
\psfrag{q}{$q$}
\psfrag{r}{$r$}
\psfrag{s}{$s$}
\psfrag{S}{$S$}
\psfrag{c}{$c$}
\psfrag{ct}{$\tilde{c}$}
\epsfxsize=.8\textwidth
\leavevmode
\epsfbox{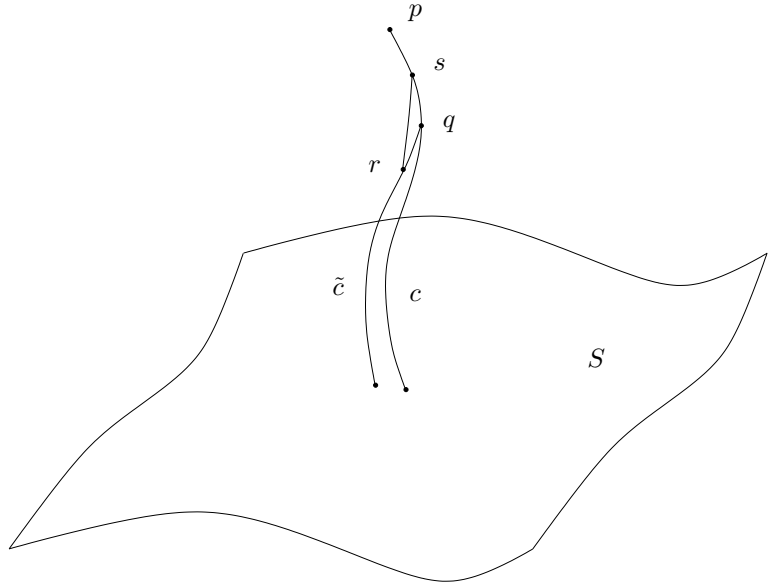}
\end{center}
\caption{Proof of Proposition~\ref{conjugate_prop}.} \label{Figure_conj}
\end{figure}

\begin{proof}
We will provide a sketch of the proof, which is similar to its analogue in Riemannian geometry. Let $q$ be a conjugate point along $c$ between $S$ and $p$. Then there exists another geodesic $\tilde{c}$, orthogonal to $S$, with the same (approximate) length, which (approximately) intersects $c$ at $q$. Let $V$ be a geodesically convex neighborhood of $q$, $r\in V$ a point along $\tilde{c}$ between $S$ and $q$, and $s\in V$ a point along $c$ between $q$ and $p$ (cf.~Figure~\ref{Figure_conj}). Then the piecewise smooth timelike curve obtained by following $\tilde{c}$ between $S$ and $r$, the unique geodesic in $V$ between $r$ and $s$, and $c$ between $s$ and $p$ connects $S$ to $p$ and has strictly bigger length than $c$ (by the Twin Paradox). This curve can be easily smoothed while retaining bigger length than $c$.
\end{proof}

\begin{Remark}
As in Riemannian geometry, a the critical values of $\exp_p$ are said to be {\bf conjugate} to $p$. Essentially the same proof as that of Proposition~\ref{conjugate_prop} shows that if $c$ is a timelike geodesic connecting $p$ to $q$ and there exists a conjugate point between $p$ and $q$ then $c$ does not maximize length (among the timelike curves connecting $p$ to $q$). An example is provided by the points $p$ and $q$ with coordinates $(0,0)$ and $(\pi + \varepsilon, 0)$ in $2$-dimensional anti-de Sitter space (cf.~Example~\ref{Hopf-Rinow}(\ref{2dAdS})). As another interesting example, consider circular motions with constant angular velocity in the Schwarzschild solution (not necessarily geodesics). From the condition that the tangent vector must be a unit timelike vector,
\[
- \left( 1 - \left(\frac{r_s}{r}\right)^{n-2} \right) \left(\frac{dt}{d\tau}\right)^2 + r^2 \left(\frac{d\varphi}{d\tau}\right)^2 = -1,
\]
one easily obtains
\[
\frac{d\tau}{dt} = \left(\left( 1 - \left(\frac{r_s}{r}\right)^{n-2} \right) - r^2 \left(\frac{d\varphi}{dt}\right)^2\right)^{\frac12},
\]
and hence the proper time along one of these motions is
\[
\Delta \tau =  \left(\left( 1 - \left(\frac{r_s}{r}\right)^{n-2} \right) - r^2 \left(\frac{d\varphi}{dt}\right)^2\right)^{\frac12} \Delta t.
\] 
Notice that $\Delta \tau$ decreases as $\left|\frac{d\varphi}{dt}\right|$ increases. A circular orbit is obtained by setting
\[
\left( \frac{d\varphi}{dt} \right)^2 = \frac{(n-2){r_s}^{n-2}}{2r^n}.
\]
Consider two points on a circular orbit more than half an orbit away. Then the (non-geodesic) circular motion in the opposite direction connecting the same points has a smaller value of $\left|\frac{d\varphi}{dt}\right|$, and hence a larger value of $\Delta \tau$. This is to be expected, as the midway point of the circular orbit is clearly conjugate to the starting point: all circular geodesics through the starting point (in different orbital planes) have the same midway point.
\end{Remark}

\begin{Prop} \label{compact_prop}
Let $(M,g)$ be a globally hyperbolic Lorentzian manifold, $S$ a Cauchy hypersurface and $p \in D^+(S)$. Then $D^+(S)\cap J^-(p)$ is compact.
\end{Prop}

\begin{proof}
Let us define a {\bf simple neighborhood} $U \subset M$ to be a geodesically convex open set diffeomorphic to an open ball whose boundary is a compact submanifold of a geodesically convex open set (therefore $\partial U$ is diffeomorphic to $S^n$ and $\overline{U}$ is compact). It is clear that simple neighborhoods form a basis for the topology of $M$. Also, it is easy to show that any open cover $\{ V_\alpha \}_{\alpha \in A}$ has a countable, locally finite refinement $\{ U_k \}_{k \in \bbN}$ by simple neighborhoods.

\begin{figure}[h!]
\begin{center}
\psfrag{p=p1}{$p=p_1$}
\psfrag{p2}{$p_2$}
\psfrag{p3}{$p_3$}
\psfrag{Un1}{$U_{k_1}$}
\psfrag{Un2}{$U_{k_2}$}
\epsfxsize=.7\textwidth
\leavevmode
\epsfbox{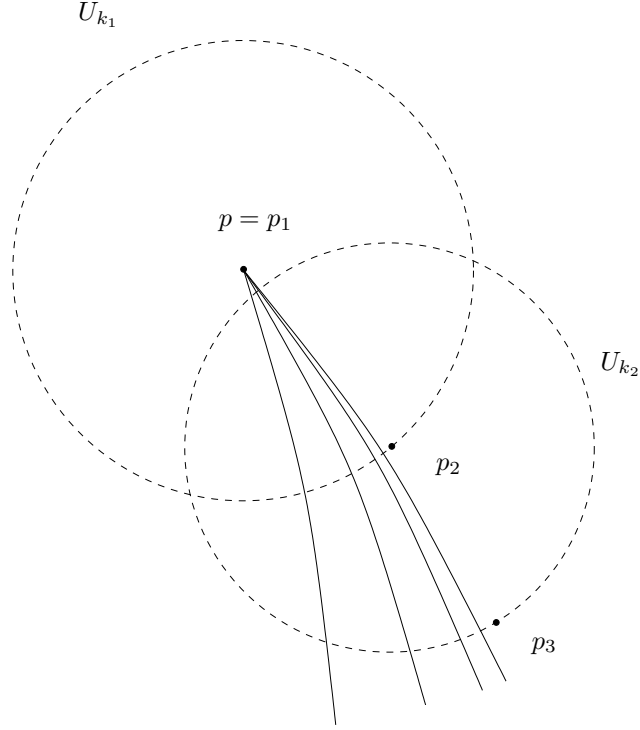}
\end{center}
\caption{Proof of Proposition~\ref{compact_prop}.} \label{Figure_comp}
\end{figure}

If $A=D^+(S)\cap J^-(p)$ were not compact then there would exist a countable, locally finite open cover $\{ U_k \}_{k \in \bbN}$ of $A$ by simple neighborhoods not admitting any finite subcover. Take $q_k \in A \cap U_k$ such that $q_k \neq q_l$ for $k \neq l$. The sequence $\{ q_k \}_{k \in \bbN}$ cannot have accumulation points, for any point in $M$ has a neighborhood intersecting only finitely many simple neighborhoods $U_k$. Consequently, each simple neighborhood $U_k$ contains only finitely many points in the sequence (as $\overline{U_k}$ is compact).

Set $p_1=p$. Since $p_1 \in A$, we have $p_1 \in U_{k_1}$ for some $k_1 \in \bbN$. Let $q_k \not\in U_{k_1}$. Since $q_k \in J^-(p_1)$, there exists a future-directed causal curve $c_k$ connecting $q_k$ to $p_1$. This curve will necessarily intersect $\partial U_{k_1}$. Let $r_{1,k}$ be an intersection point. Since $U_{k_1}$ contains only finitely many points in the sequence $\{ q_k \}_{k \in \bbN}$, there will exist infinitely many intersection points $r_{1,k}$. As $\partial U_{k_1}$ is compact, these will accumulate to some point $p_2 \in \partial U_{k_1}$.

Because $\overline{U_{k_1}}$ is contained in a geodesically convex open set, $p_2 \in J^-(p_1)$: if $\gamma_{1,k}$ is the unique causal geodesic connecting $p_1$ to $r_{1,k}$, parametrized by the global time function $t:M \to \bbR$ such that $S=t^{-1}(0)$, then the subsequence of $\{\gamma_{1,k}\}$ corresponding to a convergent subsequence of $\{r_{1,k}\}$ will converge (in the $C^\infty$ topology) to a causal geodesic $\gamma_1$ connecting $p_1$ to $p_2$. Since $t(r_{1,k}) \geq 0$, we have $t(p_2) \geq 0$, and therefore $p_2 \in A$. Since $p_2 \not\in U_{k_1}$, there must exist $k_2 \in \bbN$ such that $p_2 \in U_{k_2}$. 

Since $U_{k_2}$ contains only finitely many points in the sequence $\{ q_k \}_{k \in \bbN}$, infinitely many curves $c_k$ must intersect $\partial U_{k_2}$ to the past of $r_{1,k}$. Let $r_{2,k}$ be the intersection points. As $\partial U_{k_2}$ is compact, $\{r_{2,k}\}$ must accumulate to some point $p_3 \in \partial U_{k_2}$. Because $\overline{U_{k_2}}$ is contained in a geodesically convex open set, $p_3 \in J^-(p_2)$: if $\gamma_{2,n}$ is the unique causal geodesic connecting $r_{1,k}$ to $r_{2,k}$, parametrized by the global time function, then the subsequence of $\{\gamma_{2,k}\}$ corresponding to convergent subsequences of both $\{r_{1,k}\}$ and $\{r_{2,k}\}$ will converge to a causal geodesic connecting $p_2$ to $p_3$. Since $J^-(p_2) \subset J^-(p_1)$, $p_2 \in A$.

Iterating the procedure above, we can construct a sequence $\{p_i\}_{i\in\bbN}$ of points in $A$ satisfying $p_i \in U_{n_i}$ with $n_i \neq n_j$ if $i \neq j$, such that $p_{i}$ is connected $p_{i+1}$ by a causal geodesic $\gamma_i$. It is clear that $\gamma_i$ cannot intersect $S$, for $t(p_{i+1}) > t(p_{i+2}) \geq 0$. On the other hand, the piecewise smooth causal curve obtained by joining the curves $\gamma_i$ can easily be smoothed into a past-directed causal curve starting at $p_1$ which does not intersect $S$. Finally, such curve is inextendible: it cannot converge to any point, as $\{p_i\}_{i\in\bbN}$ cannot accumulate. But since $p_1 \in D^+(S)$, this curve would have to intersect $S$. Therefore $A$ must be compact.
\end{proof}

\begin{Cor}
Let $(M,g)$ be a globally hyperbolic Lorentzian manifold and $p,q \in M$. Then:
\begin{enumerate}[(i)]
\item
$J^+(p)$ is closed;
\item
$J^+(p) \cap J^-(q)$ is compact.
\end{enumerate}
\end{Cor}

Proposition~\ref{compact_prop} is a key ingredient in establishing the following fundamental result:

\begin{Thm} \label{max_thm}
Let $(M,g)$ be a globally hyperbolic Lorentzian manifold with Cauchy hypersurface $S$, and $p \in D^+(S)$. Then among all timelike curves connecting $p$ to $S$ there exists a timelike curve with maximal length. This curve is a timelike geodesic, orthogonal to $S$.
\end{Thm}

\begin{figure}[h!]
\begin{center}
\psfrag{p}{$p$}
\psfrag{c}{$c$}
\psfrag{g}{$\gamma$}
\psfrag{S}{$S$}
\psfrag{u=0}{$\cT=0$}
\psfrag{u=t(c)}{$\cT=\tau(c)$}
\psfrag{Ue(c)}{$U_{\varepsilon}(c)$}
\epsfxsize=.6\textwidth
\leavevmode
\epsfbox{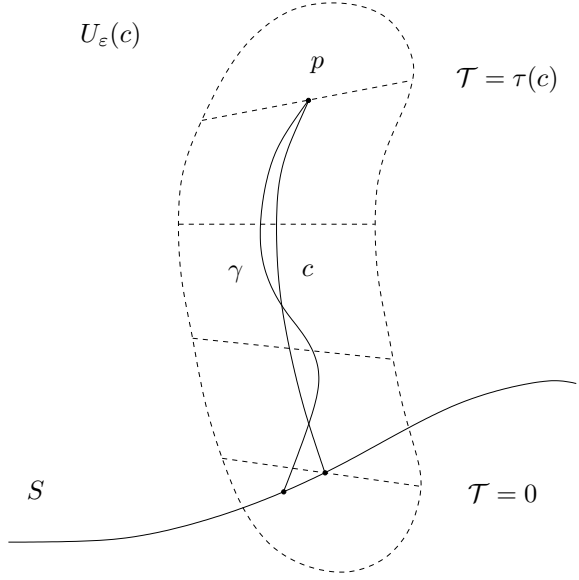}
\end{center}
\caption{Proof of Theorem~\ref{max_thm}.} \label{Figure_semi}
\end{figure}

\begin{proof}
Consider the set $T(S,p)$ of all timelike curves connecting $S$ to $p$. Since we can always use the global time function $t:M\to \bbR$ such that $S=t^{-1}(0)$ as a parameter, these curves are determined by their images, which are compact subsets of the compact set $A=D^+(S) \cap J^-(p)$. As is well known (cf.~\cite{Munkres00}), the set $C(A)$ of all compact subsets of $A$ is a compact metric space for the {\bf Hausdorff metric} $d_H$, defined as follows: if $d:M\times M \to \bbR$ is a metric yielding the topology of $M$,
\[
d_H(K,L) = \inf \{ \varepsilon > 0 \mid K \subset U_\varepsilon(L) \text{ and } L \subset U_\varepsilon(K) \},
\]
where $U_{\varepsilon}(K)$ is a $\varepsilon$-neighborhood of $K$ for the metric $d$. Therefore, the closure $C(S,p) = \overline{T(S,p)}$ is a compact subset of $C(A)$. It is not difficult to show that $C(S,p)$ can be identified with the set of {\bf continuous causal curves} connecting $S$ to $p$ (a continuous curve $c:[0,t(p)]\to M$ is said to be {\bf causal} if $c(t_2) \in J^+(c(t_1))$ whenever $t_2 > t_1$).

The length function $\tau:T(S,p) \to \bbR$ is defined by
\[
\tau(c)=\int_0^{t(p)} |\dot{c}(t)| dt.
\]
This function is {\bf upper semicontinuous}, i.e.~continuous for the topology 
\[
\mathcal{O}=\{(-\infty,a) \mid -\infty \leq a \leq +\infty \}
\]
of $\bbR$. Indeed, let $c \in T(S,p)$ be parameterized by its arclength $\cT$. For a sufficiently small $\varepsilon > 0$, the function $\cT$ can be extended to the $\varepsilon$-neighborhood $U_{\varepsilon}(c)$ in such a way that its level hypersurfaces are spacelike and orthogonal to $c$, that is, $-\grad \cT$ is timelike and coincides with $\dot{c}$ on $c$ (cf.~Figure~\ref{Figure_semi}). If $\gamma \in T(S,p)$ is in the open ball $B_\varepsilon(c) \subset C(A)$ for the Hausdorff metric $d_H$ then we can use $\cT$ as a parameter, thus obtaining
\[
d\cT(\dot{\gamma}) = 1 \Leftrightarrow \langle \dot{\gamma}, \grad \cT \rangle = 1.
\]
Therefore $\dot{\gamma}$ can be decomposed as
\[
\dot{\gamma} = \frac1{\langle \grad \cT, \grad \cT \rangle} \grad \cT + X,
\]  
where $X$ is spacelike and orthogonal to $\grad \cT$, and so
\[
|\dot{\gamma}| = \left| \frac1{\langle \grad \cT, \grad \cT \rangle} + \langle X,X \rangle \right|^\frac12.
\]
Given $\delta > 0$, we can choose $\varepsilon>0$ sufficiently small so that
\[
-\frac1{\langle \grad \cT, \grad \cT \rangle} < \left(1 + \frac{\delta}{2\tau(c)}\right)^2
\]
on the $\varepsilon$-neighborhood $U_{\varepsilon}(c)$ (as $\langle \grad \cT, \grad \cT \rangle=-1$ on $c$). We have
\[
\tau(\gamma) =  \int_0^{t(p)} \left|\frac{d \gamma}{dt} \right| \, dt = \int_0^{t(p)} |\dot{\gamma}| \frac{d\cT}{dt} \, dt = \int_{\cT(\gamma\cap S)}^{\tau(c)} |\dot{\gamma}| \, d \cT,
\]
where we have to allow for the fact that $c$ is not necessarily orthogonal to $S$, and so the initial endpoint of $\gamma$ is not necessarily at $\cT=0$ (cf.~Figure~\ref{Figure_semi}). Consequently,
\begin{align*}
\tau(\gamma) & = \int_{\cT(\gamma\cap S)}^{\tau(c)} \left| -\frac1{\langle \grad \cT, \grad \cT \rangle} - \langle X,X \rangle \right|^\frac12 \, d \cT \\
& < \int_{\cT(\gamma\cap S)}^{\tau(c)} \left(1 + \frac{\delta}{2\tau(c)}\right) \, d \cT  = \left(1 + \frac{\delta}{2\tau(c)}\right) \left(\tau(c) - \cT(\gamma\cap S)\right).
\end{align*}
Choosing $\varepsilon$ sufficiently small so that 
\[
|\cT|< \left( \frac1{\tau(c)} + \frac{2}{\delta} \right)^{-1}
\]
on $S \cap U_{\varepsilon}(c)$, we obtain $\tau(\gamma) < \tau(c) + \delta$, proving upper semicontinuity in $T(S,p)$.

As a consequence, the length function can be extended to $C(S,p)$ through
\[
\tau(c)=\lim_{\varepsilon \to 0} \sup\{ \tau(\gamma) \mid \gamma \in B_\varepsilon(c) \cap T(S,p) \}
\]
(as for $\varepsilon>0$ sufficiently small the supremum will be finite). Also, it is clear that if $c \in T(S,p)$ then the upper semicontinuity of the length forces the two definitions of $\tau(c)$ to coincide. The extension of the length function to $C(S,p)$ is trivially upper semicontinuous: given $c \in C(S,p)$ and $\delta > 0$, let $\varepsilon>0$ be such that $\tau(\gamma) < \tau(c) + \frac{\delta}2$ for any $\gamma \in B_{2\varepsilon}(c) \cap T(S,p)$. Then it is clear that $\tau(c') < \tau(c) + \delta$ for any $c' \in B_{\varepsilon}(c)$.

Finally, we notice that the compact sets of $\bbR$ for the topology $\mathcal{O}$ are sets with maximum. Therefore, the length function attains a maximum at some point $c \in C(S,p)$. All that remains to be seen is that the maximum is also attained at a smooth timelike curve $\gamma$. To do so, cover $c$ with finitely many geodesically convex neighborhoods and choose points $p_1, \ldots, p_m$ in $c$ such that $p_1 \in S$, $p_m=p$ and the portion of $c$ between $p_{i-1}$ and $p_{i}$ is contained in a geodesically convex neighborhood for all $i=2, \ldots, m$. It is clear that there exists a sequence $c_k \in T(S,p)$ such that $c_k \to c$ and $\tau(c_k) \to \tau(c)$. Let $t_i=t(p_i)$ and $p_{i,k}$ be the intersection of $c_k$ with $t^{-1}(t_i)$. Replace $c_k$ by the sectionally geodesic curve $\gamma_k$ obtained by joining $p_{i-1,k}$ to $p_{i,k}$ in the corresponding geodesically convex neighborhood. Then $\tau(\gamma_k) \geq \tau(c_k)$, and therefore $\tau(\gamma_k) \to \tau(c)$. Since each sequence $p_{i,k}$ converges to $p_i$, $\gamma_k$ converges to the sectionally geodesic curve $\gamma$ obtained by joining $p_{i-1}$ to $p_{i}$ ($i=2, \ldots, m$), and it is clear that $\tau(\gamma_k) \to \tau(\gamma)=\tau(c)$. Therefore $\gamma$ is a point of maximum for the length. Finally, we notice that $\gamma$ must be smooth at the points $p_i$, for otherwise we could increase its length by using the Twin Paradox. Therefore $\gamma$ must be a timelike geodesic. Using the Gauss Lemma, it is clear that $\gamma$ must be orthogonal to $S$, for otherwise it would be possible to increase its length.
\end{proof}

\begin{Remark}
Essentially the same proof can be used to show that if $(M,g)$ is globally hyperbolic and $p, q \in M$ with $q \in I^+(p)$ then among all timelike curves connecting $p$ to $q$ there exists a timelike curve with maximal length, which is a timelike geodesic (note that the Lorentzian manifolds in Examples~\ref{not_global_hyp} and \ref{Hopf-Rinow}(\ref{2dAdS}) are not globally hyperbolic). Thus, in a way, global hyperbolicity is the Lorentzian analogue of completeness in Riemannian geometry.
\end{Remark}

We have now all the necessary ingredients to prove the singularity theorem:

\begin{Thm} \label{sing_thm}
Let $(M,g)$ be a globally hyperbolic Lorentzian manifold satisfying the strong energy condition, and suppose that the expansion satisfies $\theta\leq\theta_0 < 0$ on a Cauchy hypersurface $S$. Then $(M,g)$ is singular.
\end{Thm}

\begin{proof}
We will show that no future-directed timelike geodesic orthogonal to $S$ can be extended to proper time greater than $\tau_0=-\frac{n}{\theta_0}$ to the future of $S$. Suppose that this was not so. Then there would exist a future-directed timelike geodesic $c$ orthogonal to $S$ defined in an interval $[0,\tau_0+\varepsilon]$ for some $\varepsilon>0$. Let $p=c(\tau_0+\varepsilon)$. According to Theorem~\ref{max_thm}, there would exist a timelike geodesic $\gamma$ with maximal length connecting $S$ to $p$, orthogonal to $S$. Because $\tau(c)=\tau_0+\varepsilon$, we would necessarily have $\tau(\gamma)\geq\tau_0+\varepsilon$. Proposition~\ref{exist_conj_prop} guarantees that $\gamma$ would develop a conjugate point at a distance of at most $\tau_0$ to the future of $S$, and Proposition~\ref{conjugate_prop} states that $\gamma$ would cease to be maximizing beyond this point. Therefore we arrive at a contradiction.
\end{proof}

\begin{Remark}
It should be clear that $(M,g)$ is singular if the condition $\theta\leq\theta_0 < 0$ on a Cauchy hypersurface $S$ is replaced by the condition $\theta\geq\theta_0 > 0$ on $S$. In this case, no {\bf past-directed} timelike geodesic orthogonal to $S$ can be extended to proper time greater than $\tau_0=\frac{n}{\theta_0}$ to the {\bf past} of $S$.
\end{Remark}

\begin{Remark}
The proof of Theorem~\ref{sing_thm} does not hold in Riemannian geometry, where geodesics are length {\bf minimizing} curves. Using similar ideas, one can prove the following: if $(M,g)$ is a complete Riemannian manifold whose Ricci curvature satisfies $Ric \geq \varepsilon g$ for some $\varepsilon > 0$ then $M$ is compact.
\end{Remark}

\begin{figure}[h!]
\begin{center}
\psfrag{t}{$x^0$}
\psfrag{x}{$x^1$}
\psfrag{y}{$x^n$}
\psfrag{incomplete}{incomplete timelike geodesic}
\psfrag{tau=constant}{$\{\tau=\text{constant}\}$}
\epsfxsize=.7\textwidth
\leavevmode
\epsfbox{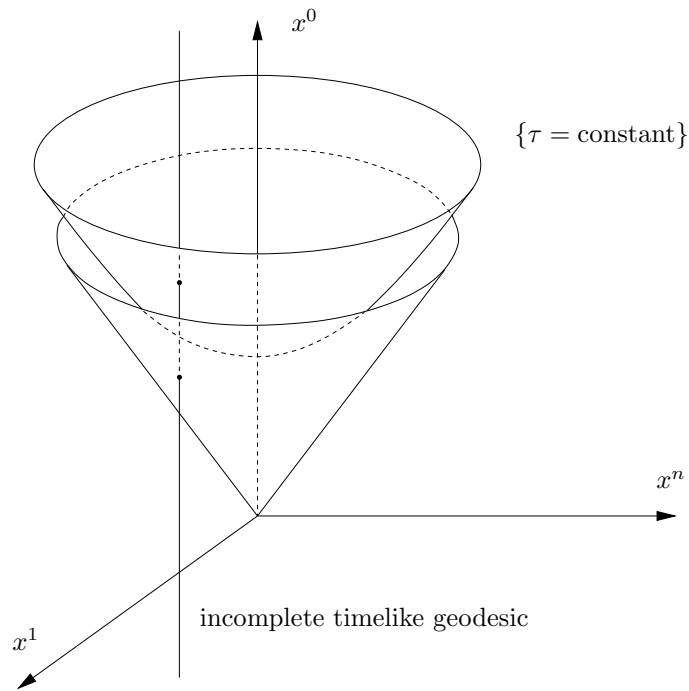}
\end{center}
\caption{Milne universe.} \label{Milne}
\end{figure}

\begin{Example} \hspace{1cm}
\begin{enumerate}
\item
Theorem~\ref{sing_thm} does not apply to Minkowski, de Sitter or anti-de Sitter spaces: the first does not contain a Cauchy hypersurface with expansion bounded away from zero, the second does not satisfy the strong energy condition and the third is not globally hyperbolic. It does apply to the {\bf Milne universe}, which is the open set $U$ of Minkowski space defined by
\[
x^0 > 0 \quad \text{ and } \quad - \left( x^0 \right)^2 + \left( x^1 \right)^2 + \ldots + \left( x^n \right)^2 < 0
\]
(cf.~Figure~\ref{Milne}). It is easily seen that the function
\[
\tau\left( x^0, x^1, \ldots, x^n \right) = \sqrt{\left( x^0 \right)^2 - \left( x^1 \right)^2 - \ldots - \left( x^n \right)^2}
\]
is a global time function in $U$ whose time slices have domain of dependence $U$ and constant positive expansion $n \tau$ (they are, in fact, isometric to $n$-dimensional hyperbolic spaces of radii $\tau$). Here the geodesic incompleteness guaranteed by Theorem~\ref{sing_thm} is not due to curvature singularities.
\item
The Friedmann-Lema\^\i tre-Robertson-Walker models are globally hyperbolic, and for $\alpha > 0$ satisfy the strong energy condition (as $\rho > 0$, cf.~Example~\ref{spacetimes}). Furthermore, one easily checks that the expansion of the time slices is
\[
\theta = \frac{n\dot{a}}{a}.
\]
Assume that the model is expanding at time $t_0$. Then $\theta = \theta_0 = \frac{n\dot{a}(t_0)}{a(t_0)}$ on the Cauchy hypersurface $S = \{ t=t_0 \}$, and hence Theorem~\ref{sing_thm} guarantees that this model is singular to the past of $S$ (i.e.~there exists a Big Bang). Furthermore, Theorem~\ref{sing_thm} implies that this singularity is generic: any sufficiently small perturbation of an expanding Friedmann-Lema\^\i tre-Robertson-Walker model satisfying the strong energy condition will also be singular. Loosely speaking, any expanding universe must have begun at a Big Bang.
\item
The region $\{ r < r_s \}$ of the Schwarzschild solution is globally hyperbolic, and satisfies the strong energy condition (as $Ric=0$). The metric can be written is this region as
\[
g = - d\tau \otimes d\tau + \left( \left(\frac{r_s}r \right)^{n-2} - 1 \right) dt \otimes dt + r^2 h,
\]
where $h$ is the round metric in $S^{n-1}$ and
\[
\tau = \int_r^{r_s} \left( \left(\frac{r_s}u \right)^{n-2} - 1 \right)^{-\frac12} du.
\]
Therefore the inside of a Schwarzschild black hole can be pictured as a cylinder $\bbR \times S^{n-1}$ whose shape is evolving in time: as $r \to 0$, the sphere $S^{n-1}$ contracts to a singularity, with the $t$-direction expanding. The expansion of the Cauchy hypersurface $S = \{ \tau=\tau_0 \} = \{ r=r_0 \}$ can be computed to be
\[
\theta = \left( \left(\frac{r_s}{r_0} \right)^{n-2} - 1 \right)^{-\frac12}\left( \frac{n-1}{r_0} - \frac{n}{2r_0} \left(\frac{r_s}{r_0} \right)^{n-2} \right).
\]
Therefore $\theta = \theta_0 < 0$ for $r_0$ sufficiently small, and hence Theorem~\ref{sing_thm} guarantees that the Schwarzschild solution is singular to the future of $S$. Moreover, Theorem~\ref{sing_thm} implies that this singularity is generic: any sufficiently small perturbation of the Schwarzschild solution satisfying the strong energy condition will also be singular. Loosely speaking, once the collapse has advanced long enough nothing can prevent the formation of a singularity.
\end{enumerate}
\end{Example}

%
%

\bibliographystyle{amsalpha} 

\end{document}